\documentclass[12pt,a4paper]{article}

\usepackage{amsthm,amsmath,natbib,amssymb,amsfonts,bm}
\usepackage{graphicx}
\usepackage{placeins}
\usepackage{subcaption}

\usepackage{sectsty} 
\sectionfont{\fontsize{12}{16}\selectfont} 
\subsectionfont{\fontsize{12}{16}\selectfont} 

\def\hang{\hangindent\parindent}
\def\rf{\par\noindent\hang}

\newtheorem{theorem}{Theorem}

\newtheorem*{theorem*}{Theorem}

\theoremstyle{definition}

\theoremstyle{remark}

\newcommand{\bs}{\boldsymbol}

\DeclareMathAlphabet{\mathpzc}{OT1}{pzc}{m}{it}

\newcommand{\bbeta}{\bs{\beta}}

\def\hang{\hangindent\parindent}
\def\rf{\par\noindent\hang}

\begin{document}

\baselineskip=20pt

\begin{center}
{\bf \Large Two sources of poor coverage of confidence intervals after model selection}
\end{center}


\bigskip

\begin{center}
{\bf \large Paul Kabaila$^*$ and Rheanna Mainzer}
\end{center}

\medskip

\begin{center}
{\sl Department of Mathematics and Statistics, La Trobe University, Australia}
\end{center}

\vspace{1cm}

\noindent \textbf{ABSTRACT}


 \medskip
 

\noindent We compare the following two sources of poor coverage of post-model-selection 
confidence intervals: the preliminary data-based model selection sometimes chooses the wrong model and the data used to choose the model is re-used for the 
construction of the confidence interval.

\bigskip
\bigskip


\vspace{5.5cm}

\noindent * Corresponding author. Department of Mathematics and Statistics,
La Trobe University, Victoria 3086, Australia. Tel.: +61 3 9479 2594; fax +61 3 9479 2466.
{\sl E-mail address:} P.Kabaila@latrobe.edu.au.

\newpage


\noindent \textbf{1. Introduction}

\medskip

In applied statistics a model is often chosen using a data-based model selection procedure.  
Such procedures include hypothesis testing, minimizing a criteron such as AIC or BIC, and maximizing a criterion such as $R^2_{adj}$  (Sheather, 2008).
Inferences made using a model that has been selected using data-based model selection can be misleading (see e.g. Chatfield 1995).  Suppose that the inference of interest is a confidence interval for a parameter $\theta$, with minimum coverage $1-\alpha$.  
Typically, a confidence interval for $\theta$ that is constructed using a model chosen by some form of data-based model selection procedure will have poor coverage properties (Kabaila, 2009, Leeb and P\"otscher, 2005).
This poor coverage probability can be attributed to two different sources.  
The first source is that the preliminary model selection procedure sometimes chooses the wrong model.  When this happens, the incorrect model is treated as if it were the true model during the construction of the confidence interval.
The second source is that the data used to choose the model is re-used for the 
construction of the confidence interval, without due acknowledgment.
It is important to know the relative importance of these two sources
of poor coverage.

Miller (2002, Ch.6) considers the problems that arise for parameter estimation 
due to the re-use of the data used to choose a model. 
Hurvich and Tsai (1990) examine the coverage probability of confidence regions for the regression parameter vector in linear regression models after data-based model selection.  They acknowledge the problem of re-using the same data and propose a possible solution:
\begin{quote}
One possible solution 
...
is to actually perform model selection and inference on separate parts of the data. 
... This procedure may seem wasteful, but we note that linear regression model 
selection can now be performed on quite small samples. ... Thus the available data set
need not be split in half; instead, a small part may be used for the model selection stage and the remainder may be used for inference. 
\end{quote}
Thus if the second source of poor coverage (re-use of the data) dominates the first source (possible choice of the wrong model) we could split the data into two parts: the first part could be used to choose the model and the second to construct a confidence interval. By allocating most of the data to the second part, we might hope to avoid constructing confidence intervals that are too wide.

If, on the other hand, the first source of poor coverage dominates the second source, this data splitting strategy would not correct the defect of poor coverage probability.  Instead we would need to consider other methods to improve the coverage probability. 
We expect that problems caused by possibly choosing a wrong model would be 
ameliorated by using an appropriate data-weighted average of models in place of putting all of the weight on a model that may be wrong.
In other words, if the first source of poor coverage dominates the second source then some form of frequentist model averaged confidence interval (Buckland, Burnham and Augustin, 1997, Hjort and Claeskens 2003,  Fletcher and Turek, 2011, Turek and Fletcher, 2012, Efron, 2014)
would be expected to greatly improve the coverage probability.

In this paper we examine the extent to which each of these two sources contributes to the difference between true and nominal coverage probability.  
The basic argument presented in Section 2 gives us the means to identify $D_{wm}$,
the deficit in coverage probability due to the possible choice of the \textbf{w}rong \textbf{m}odel, 
and $D_{rd}$, the deficit in coverage probability due to \textbf{r}e-use of the \textbf{d}ata.
The coverage probability of the confidence interval, constructed after a preliminary data-based model selection, is given by $1 - \alpha - D_{wm} - D_{rd}$.

In Section 3, we apply this basic argument to
the simple scenario that preliminary data-based
choice is between two nested linear regression models.  The simpler model is obtained from the full model by setting $\tau$, a given linear combination of the regression parameters, to 0.  The parameter of interest $\theta$ is a distinct given linear combination of the regression parameters.  Let $\widehat{\theta}$ and $\widehat{\tau}$ denote the least squares estimators of $\theta$ and $\tau$, respectively.  
We consider preliminary model selection using a t test of the null hypothesis $\tau = 0$ against the alternative hypothesis $\tau \ne 0$. 
Let $n$ denote the dimension of the response vector and $p$ denote the dimension of the regression parameter vector. 
Define the unknown parameter $\gamma = \tau / (\text{var}(\widehat{\tau}))^{1/2}$ and the known correlation $\rho = \text{corr}(\widehat{\theta}, \widehat{\tau})$.

\begin{figure}[h]
	\centering
\includegraphics[scale = 0.95]{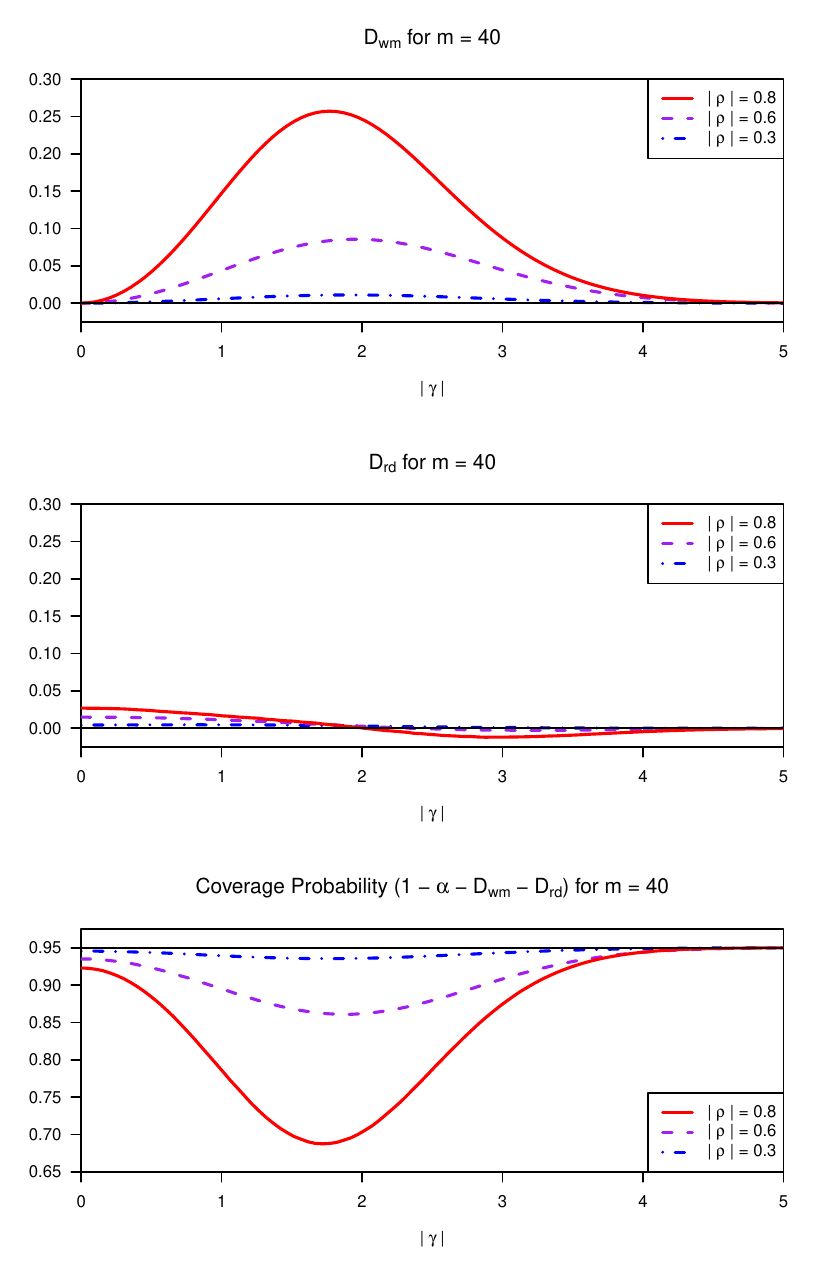}
	\caption{This figure is for the confidence interval, with nominal coverage 0.95, constructed after a preliminary t test with size 0.1, when $n - p = 40$
		and $|\rho| \in \{0.3, 0.6, 0.8 \}$. 
		The top panel shows graphs of $D_{wm}$, the deficit in coverage due to the choice of the wrong model, as a function of $|\gamma|$.
		The middle panel shows graphs of $D_{rd}$, the deficit in coverage due to the re-use of the data, as a function of $|\gamma|$.
		The bottom panel shows graphs of the coverage probability of this confidence interval,
		given by $1 - \alpha - D_{wm} - D_{rd}$, as a function of $|\gamma|$.}
	\label{Fig1}
\end{figure}

Figure \ref{Fig1} is for the confidence interval, with nominal coverage 0.95, constructed after a preliminary t test with size 0.1, when $n - p = 40$
and $|\rho| \in \{0.3, 0.6, 0.8 \}$. 
The top panel shows graphs of $D_{wm}$, the deficit in coverage due to the choice of the wrong model, as a function of $|\gamma|$.
The middle panel shows graphs of $D_{rd}$, the deficit in coverage due to the re-use of the data, as a function of $|\gamma|$.
The bottom panel shows graphs of the coverage probability of this 
post-model-selection confidence interval,
given by $1 - \alpha - D_{wm} - D_{rd}$, as a function of $|\gamma|$.
We do not include the value  $\rho = 0$ in these 
graphs since, to a very good approximation, $D_{wm} = 0$, $D_{rd} = 0$ and  the 

\FloatBarrier

\noindent coverage probability is 0.95 for all $|\gamma|$. 
This figure shows that, for size 0.1 of the preliminary t test, $n - p = 40$
 and $|\rho| \in \{0.6, 0.8\}$,
$D_{wm}$, the deficit in coverage due to the choice of the wrong model, dominates
$D_{rd}$, the deficit in coverage due to the re-use of the data.
An extensive numerical 
investigation of $D_{wm}$ and $D_{rd}$, covering a wide range of values of $n-p$, $\rho$, $1 - \alpha$ 
and size of the preliminary t test, is reported in the Supplementary Material and summarised in Section 4.

\bigskip


\noindent \textbf{2. The basic argument}

\medskip

Suppose that we observe a random $n$-vector of responses  $\bs{Y}$ for a general 
regression model, with given values of the explanatory variables. We denote this 
``full model'' by ${\cal M}$. Suppose that the parameter of interest is the
scalar $\theta$. Also suppose that a given preliminary data-based model selection 
is used to select between the full model and a specified set of submodels.
Assuming that the selected model had been given to us \textit{a priori} as the true
model, we construct the usual confidence interval for $\theta$, with nominal coverage
$1 - \alpha$. Let $K$ denote the resulting post-model-selection confidence interval.

To identify the deficit in coverage probability due to the possible choice of the wrong model and the deficit in coverage probability due to re-use of the data, we consider the
following hypothetical situation. Suppose that a genie provides an observation of 
$\bs{Y}^*$ the random $n$-vector of responses for a replicate (i.e. an independent realization) of the original experiment that gave rise to $\bs{Y}$, with the same 
values of the explanatory variables. Also suppose that the same preliminary model selection based on $\bs{Y}$ is used. Assuming that the selected model had been given to us \textit{a priori} as the true
model, we construct the usual confidence interval for $\theta$, with nominal coverage
$1 - \alpha$ using the genie's data $\bs{Y}^*$. Let $K^*$ denote the resulting confidence interval.

Observe that $K$ and $K^*$ are influenced by exactly the same, possibly wrong, choice of model based on $\bs{Y}$. However, unlike $K$, the confidence interval
$K^*$ is not influenced by the re-use of the data for its construction.
 Because of this difference
between $K$ and $K^*$, we are able to separate out these two sources of poor coverage as follows. The coverage probability of $K$ is 
\begin{align*}
P\big(\theta \in K \big) 
%
&= (1 - \alpha) - \Big( (1 - \alpha) - P \big( \theta \in K^* \big) \Big) - \Big( P \big( \theta \in K^* \big) - P \big(\theta \in K \big) \Big)  \\[1.2ex]
%
&= (1 - \alpha) - D_{wm} - D_{rd},
\end{align*}
where $D_{wm} = (1 - \alpha) - P\big(\theta \in K^* \big)$ is the deficit in coverage probability due solely to possibly choosing the wrong model and $D_{rd} = P \big( \theta \in K^* \big) - P \big( \theta \in K \big)$ is the deficit in coverage probability due solely to re-using the data. 

\bigskip

\noindent \textbf{3. The scenario of two nested linear regression models}

\medskip


Consider the multiple linear regression model 
\begin{equation}
\label{full_model}
\bs{Y} = \bs{X} \bs{\beta} + \bs{\varepsilon}
\end{equation}
where $\bs{Y}$ is a random $n$-vector of responses, $\bs{X}$ is a known $n \times p$ matrix with linearly independent columns, $\bs{\beta}$ is an unknown parameter $p$-vector and $\bs{\varepsilon} \sim N(\bs{0}, \sigma^2 \bs{I}_n)$, where $\sigma^2$ is an unknown positive parameter.   Suppose that $n > p$. 
 We refer to \eqref{full_model} as the full model and denote it by ${\cal M}$.

Suppose that the quantity of interest is $\theta = \bs{a}^{\top} \bs{\beta}$, where $\bs{a}$ is a known $p$-vector ($\bs{a} \not = \bs{0}$).  Also suppose that the inference of interest is a $1 - \alpha$ confidence interval for $\theta$.  Define the parameter $\tau = \bs{c}^{\top} \bs{\beta} - r$, where $\bs{c}$ is a known $p$-vector, the number $r$ is specified and the vectors $\bs{a}$ and $\bs{c}$ are linearly independent.   
Let $\widehat{\bs{\beta}}$ denote the least squares estimator of $\bs{\beta}$.  Then $\widehat{\theta} = \bs{a}^{\top} \, \widehat{\bs{\beta}}$ is the least squares estimator of $\theta$ and $\widehat{\tau} = \bs{c}^{\top} \, \widehat{\bs{\beta}} - r$ is the least squares estimator of $\tau$.  Let $m = n - p$ and 
$\widehat{\sigma}^2 
= (\bs{Y} - \bs{X} \widehat{\bs{\beta}})^{\top} (\bs{Y} - \bs{X} \widehat{\bs{\beta}}) \big / m$.  
Also let $v_{\theta} = \text{var}(\widehat{\theta}) / \sigma^2$, $v_{\tau} = \text{var}(\widehat{\tau})  / \sigma^2$, $\rho = \text{corr}(\widehat{\theta}, \widehat{\tau})$ and $\gamma = \tau \big / (\sigma \, v_{\tau}^{1/2})$

We consider preliminary model selection using a t test of the null hypothesis $H_0: \tau = 0$ against the alternative hypothesis $H_A: \tau \not= 0$.  The test statistic is 
\begin{equation*}
T = \frac{\widehat{\tau}}{\widehat{\sigma} \,v_{\tau}^{1/2}}.
\end{equation*}
Define the quantile $t_{m, a}$ by $P(T \leq t_{m, a}) = a$ for $T \sim t_m$.  
Suppose that we reject $H_0$ if $|T| \geq t_{m, \, 1 - \widetilde{\alpha}/2}$; otherwise we accept $H_0$. This test has size $\widetilde{\alpha}$.
Let ${\cal M}_{\tau=0}$ denote the model ${\cal M}$ for $\tau = 0$. 
If we reject $H_0$ then we choose the model ${\cal M}_{\tau=0}$;
otherwise we choose the model ${\cal M}$.

\bigskip

\noindent \textbf{3.1 \ The confidence interval $\bs{K}$}

\medskip

 The usual $1 - \alpha$ confidence interval for $\theta$ based on the full model ${\cal M}$ is 
\begin{equation*}
I = \big[ \widehat{\theta} \pm t_{m, \, 1 - \alpha/2} \, v_{\theta}^{1/2} \, \widehat{\sigma} \big],
\end{equation*}
where $[a \pm b] = [a - b, a + b]$ ($b \ge 0$). 
Under the model ${\cal M}_{\tau=0}$ the usual $1 - \alpha$ confidence interval for $\theta$ is 
\begin{align*}
J = \left[ \widehat{\theta} - \rho \, v_{\theta}^{1/2} \, (\widehat{\tau}/v_{\tau}^{1/2}) \pm t_{m+1, 1 - \alpha/2} \; v_{\theta}^{1/2} \, \left( 1 - \rho^2 \right)^{1/2} \, \left( \frac{m \, \widehat{\sigma}^2 + (\widehat{\tau}^2 / v_{\tau})}{m + 1} \right)^{1/2} \, \right].
\end{align*}

Let $K$ be the confidence interval, with nominal coverage
$1 - \alpha$, for $\theta$ constructed as follows.
We use the data $\bs{Y}$ to choose the model.
We then construct the
confidence interval for $\theta$ using the same data $\bs{Y}$
based on the assumption that this model had been given to us
\textit{a priori} as the true model. Therefore $K$ is the 
 post-model-selection confidence interval, with nominal coverage
$1 - \alpha$, for $\theta$.  In other words, 
\begin{equation*}
K = 
\begin{cases} 
I &\text{ if } \ |T| \geq t_{m, 1 - \widetilde{\alpha}/2} \\ 
J &\text { if } \ |T| < t_{m, 1 - \widetilde{\alpha}/2}.
\end{cases}
\end{equation*}

\bigskip

\noindent \textbf{3.2 \ The confidence interval $\bs{K}^*$}


\medskip

Now suppose that a genie provides the following
replicate (i.e. an independent realization) of the original experiment that gave rise to $\bs{Y}$, with the same 
values of the explanatory variables:
\begin{equation*}
\bs{Y}^* = \bs{X} \bs{\beta} + \bs{\varepsilon}^*,
\end{equation*}
where $\bs{Y}^*$ is a random $n$-vector of responses and $\bs{\varepsilon}^* \sim N(\bs{0}, \sigma^2 \bs{I}_n)$.  
Note that $\bs{\varepsilon}$ and $\bs{\varepsilon}^*$ are independent.
Let $\widehat{\bbeta}^*$, $\widehat{\theta}^*$ and $\widehat{\tau}^*$ denote the least squares estimator of $\bs{\beta}$, $\theta$ and $\tau$, respectively, based on the genie's data $\bs{Y}^*$.  Let 
$(\widehat{\sigma}^*)^2 
= (\bs{Y}^* - \bs{X} \widehat{\bs{\beta}}^*)^{\top}(\bs{Y}^* - \bs{X} \widehat{\bs{\beta}}^*)/m$.

The usual $1 - \alpha$ confidence interval for $\theta$ based on ${\cal M}$, and using only the genie's data $\bs{Y}^*$, is 
\begin{equation*}
I^* = \big[ \widehat{\theta}^* \pm t_{m, \, 1 - \alpha/2} \, v_{\theta}^{1/2} \, \widehat{\sigma}^* \big]. 
\end{equation*}
The usual $1 - \alpha$ confidence interval for $\theta$ based on ${\cal M}_{\tau = 0}$, and using only the genie's data $\bs{Y}^*$, is
\begin{equation*}
J^* = \left[ \widehat{\theta}^* -  \rho \, v_{\theta}^{1/2} \, (\widehat{\tau}^*/v_{\tau}^{1/2}) \pm t_{m+1, 1 - \alpha/2} \, v_{\theta}^{1/2} \, \left( 1 - \rho^2 \right)^{1/2} \, \left( \frac{m \, (\widehat{\sigma}^*)^2 + ((\widehat{\tau}^*)^2 / v_{\tau})}{m + 1} \right)^{1/2} \, \right].
\end{equation*}

Let $K^*$ be the confidence interval, with nominal coverage
$1 - \alpha$, for $\theta$ constructed as follows.
We use the original data $\bs{Y}$ to choose the model
i.e. we use the test statistic $T$.
We then construct the
confidence interval for $\theta$ using the genie's data $\bs{Y}^*$
based on the assumption that this model had been given to us
\textit{a priori} as the true model. In other words, 
\begin{equation*}
K^* = \begin{cases} 
I^* &\text{ if } \ |T| \geq t_{m, \, 1 - \widetilde{\alpha}/2} \\
J^* &\text{ if } \ |T| < t_{m, \, 1 - \widetilde{\alpha}/2}. 
\end{cases}
\end{equation*}

\bigskip

\noindent \textbf{3.3 \ Computationally convenient formulas for the deficits $\boldsymbol{D_{wm}}$ and $\boldsymbol{D_{rd}}$}


\medskip

By the law of total probability,
\begin{align*}
P(\theta \in K) &= P(\theta \in J, \, |T| < t_{m, \, 1 - \widetilde{\alpha}/2}) + P (\theta \in I, \, |T| \geq t_{m, \, 1 - \widetilde{\alpha}/2}) \\[1.2ex]
&= P(\theta \in J, \, |T| < t_{m, \, 1 - \widetilde{\alpha}/2}) + P(\theta \in I) - P(\theta \in I, \, |T| < t_{m, \, 1 - \widetilde{\alpha}/2}) \\[1.2ex]
&= (1 - \alpha) + P(\theta \in J, \, |T| < t_{m, \, 1 - \widetilde{\alpha}/2}) - P(\theta \in I, \, |T| < t_{m, \, 1 - \widetilde{\alpha}/2}),
\end{align*}
since $P(\theta \in I) = 1 - \alpha$.
By the law of total probability, and the independence between the original data $\bs{Y}$ and the genie's data $\bs{Y}^*$,
\begin{align*}
P(\theta \in K^*) &= P (\theta \in J^*, \, |T| < t_{m, \, 1 - \widetilde{\alpha}/2}) + P(\theta \in I^*, \, |T| \geq t_{m, \, 1 - \widetilde{\alpha}/2}) \\[1.2ex] 
&= P(\theta \in J^*) \, P(|T| < t_{m, \, 1 - \widetilde{\alpha}/2}) + P(\theta \in I^*) \, P(|T| \geq t_{m, \, 1 - \widetilde{\alpha}/2}) \\[1.2ex]
&= P(\theta \in J) \, P(|T| < t_{m, \, 1 - \widetilde{\alpha}/2}) + P(\theta \in I) \left( 1 - P(|T| < t_{m, \, 1 - \widetilde{\alpha}/2}) \right) \\[1.2ex]
&= P(\theta \in I) + P(|T| < t_{m, \, 1 - \widetilde{\alpha}/2}) \big( P(\theta \in J) - P (\theta \in I) \big) \\[1.2ex]
&= (1 - \alpha) +  P(|T| < t_{m, \, 1 - \widetilde{\alpha}/2}) \big( P(\theta \in J) - (1 - \alpha) \big).
\end{align*}

Therefore, the deficit in coverage probability due solely to possibly choosing the wrong model using the original data $\bs{Y}$ is, by definition,
\begin{align*}
D_{wm} &= (1 - \alpha) - P(\theta \in K^*) \\[1.2ex]
&= P(|T| < t_{m, \, 1 - \widetilde{\alpha}/2}) \big( (1 - \alpha) - P(\theta \in J) \big),
\end{align*}
and the deficit in coverage probability due solely to re-using the original data $\bs{Y}$ is, by definition,
\begin{align*}
D_{rd} &= P(\theta \in K^*) - P(\theta \in K) \\[1.2ex]
&= P\big(\theta \in I, |T| < t_{m, \, 1 - \widetilde{\alpha}/2}\big) - P\big(\theta \in J, \, |T| < t_{m, \, 1 - \widetilde{\alpha}/2}\big) \\[1.2ex]
& \indent + P\big(|T| < t_{m, \, 1 - \widetilde{\alpha}/2}\big) \left( P\big(\theta \in J\big) - (1 - \alpha) \right).
\end{align*}

Let $W = \widehat{\sigma} / \sigma$ and note that $W$ has the same distribution as $(R/m)^{1/2}$, where $R \sim \chi^2_m$.   
Also let $G = (\widehat{\theta} - \theta) \big/ (\sigma \, v_{\theta}^{1/2})$
and $H = \widehat{\tau} \big/ (\sigma \, v_{\tau}^{1/2})$. Note that $W$ and $(G,H)$ 
are independent, with
\begin{equation}
\label{(G,H)dist}
\begin{bmatrix}
\, G \, \\ \, H \,
\end{bmatrix}
\sim 
N \left( \begin{bmatrix} 0 \, \\ \gamma \, \end{bmatrix}, \, \begin{bmatrix} 1 & \rho \\ 
\rho  & 1 \end{bmatrix} \right).
\end{equation}
Thus $T = H / W$ has a noncentral t distribution with degrees of freedom $m$ and noncentrality parameter $\gamma$.
 Therefore we can easily evaluate $P\big(|T| < t_{m, \, 1 - \widetilde{\alpha}/2}\big)$. We now show how the other probabilities that appear in the 
expressions for $D_{wm}$ and $D_{rd}$ are readily evaluated.
Let $\phi$ 
and $\Phi$ denote the $N(0,1)$ pdf and cdf, respectively.
Also
let $f_W$ denote the pdf of $W$.  
Define $\Psi(x, u) = \Phi\big(g_2(x, u)\big) - \Phi\big(g_1(x, u)\big)$, where
\begin{align*}
g_1(x, u) =  - t_{m+1, 1-\alpha/2} \;  x  \, \left(\dfrac{m + u^2}{m + 1} \right)^{1/2} + \dfrac{ \rho \, \gamma}{ \left(1 - \rho^2 \right)^{1/2}}
\end{align*}
and
\begin{align*}
g_2(x, u) = t_{m+1, 1-\alpha/2} \; x \, \left( \dfrac{m + u^2}{m + 1} \right)^{1/2} + \dfrac{\rho \, \gamma}{ \left(1 - \rho^2 \right)^{1/2}}.
\end{align*}
The following theorem gives computationally convenient formulas for the probabilities that we need to evaluate to find $D_{wm}$ and $D_{rd}$.

\begin{theorem}
\label{thm_expressions}
\begin{equation*}
P(\theta \in J ) = \int_0^{\infty} \int_{-\infty}^{\infty} \Psi(w, h/w) \, \phi(h - \gamma) \, dh \, f_W(w) \, dw,
\end{equation*}
$P\big(\theta \in J, |T| < t_{m, \, 1 - \widetilde{\alpha}/2}\big)$ is equal to
\begin{equation*}
t_{m, 1 - \widetilde{\alpha}/2} \,
\int_0^{\infty} \int_{-1}^1 \Psi \left(w, t_{m, 1 - \widetilde{\alpha}/2} \, y \right)  \, 
		 \phi(t_{m, 1 - \widetilde{\alpha}/2} \, w \, y - \gamma) \, dy \,  w  \, f_W(w)  \, dw,
\end{equation*}
and $P\big(\theta \in I, |T| < t_{m, \, 1 - \widetilde{\alpha}/2}\big)$ is equal to
\begin{equation*}
P \left( -t_{m, \, 1 - \alpha/2} \leq \dfrac{G}{W} \leq t_{m, \, 1 - \alpha/2}, \, -t_{m, \, 1- \widetilde{\alpha}/2} \leq \dfrac{H}{W} \leq t_{m, \, 1 - \widetilde{\alpha}/2} \right).
\end{equation*}
This probability is readily evaluated using the fact that $\big(G/W, \, H/W \big)$ has a bivariate noncentral
t distribution (Kshirsagar definition).

\end{theorem}

Theorem \ref{thm_expressions} is proved in the appendix.  It follows from Theorem \ref{thm_expressions} that the two deficits $D_{wm}$ and $D_{rd}$, and the coverage probability $P(\theta \in K)$ depend  only on one unknown parameter $\gamma$, and  four known quantities $1-\alpha$, $\widetilde{\alpha}$, 
$\rho$ and $m$.
Using Theorem \ref{thm_expressions}, it is simple to calculate $D_{wm}$, $D_{rd}$ and $P(\theta \in K)$. Our computer programs are written in \texttt{R}.  We use the \texttt{pmvt} function in the \texttt{R} package \texttt{mvtnorm} to evaluate the bivariate noncentral t distribution  (Kshirsagar definition).

The following theorem reduces the number of combinations of $\rho$ and $\gamma$ that we need to consider in order to understand the properties of the $D_{wm}$, $D_{rd}$ and $P(\theta \in K)$ functions.

\begin{theorem}
\label{thm_evenness}
For any given $\rho$, the deficits $D_{wm}$, $D_{rd}$ and the coverage probability $P(\theta \in K)$ are even functions of $\gamma$.  For any given $\gamma$, the deficits $D_{wm}$, $D_{rd}$ and the coverage probability $P(\theta \in K)$ are even functions of $\rho$.
\end{theorem}

Theorem \ref{thm_evenness} will be true if each of the probabilities appearing in the expressions for $D_{wm}$ and $D_{rd}$ are, for any given $\rho$, even functions of $\gamma$
and, for any given $\gamma$, even functions of $\rho$.
Appendix C of Kabaila and Giri (2009) provides a proof that the coverage probability $P(\theta \in K)$ is, for any given $\rho$, an even function of $\gamma$
and, for any given $\gamma$, an even function of $\rho$.
It follows directly from this proof that $P(\theta \in J, |T| < t_{m, 1 - \widetilde{\alpha}/2})$ and $P(\theta \in I, |T| < t_{m, 1 - \widetilde{\alpha}/2})$ are, for any given $\rho$, even functions of $\gamma$
and, for any given $\gamma$, even functions of $\rho$.
 Using the same arguments as those given in Appendix C of Kabaila and Giri (2009), it is straightfoward to show that $P(\theta \in J)$ and $P(|T| < t_{m, 1 - \widetilde{\alpha}/2})$ are, for any given $\rho$, even functions of $\gamma$
 and, for any given $\gamma$, even functions of $\rho$.
 
 \bigskip
 
 \noindent \textbf{4. Extensive numerical investigation of the deficits $\boldsymbol{D_{wm}}$ and $\boldsymbol{D_{rd}}$}
 

\medskip

In this section we summarise the extensive numerical investigation of $D_{wm}$, the deficit in coverage due to the choice of the wrong model, and $D_{rd}$, the deficit in coverage due to the re-use of the data, reported in the Supplementary Material.
Primarily our interest lies in the cases where the coverage probability of $K$, minimized over $\gamma$, i.e. $\min_{\gamma} P(\theta \in K)$, is substantially below $1 - \alpha$.  In the Supplementary Material we compute $\min_{\gamma}P(\theta \in K)$ for $m \in \{1, 2, 5, 10, 40, 100\}$, $|\rho| \in \{0, 0.3, 0.6, 0.8 \}$, $1-\alpha \in \{0.9, 0.95, 0.98\}$ and $\widetilde{\alpha} \in \{0.02, 0.05, 0.1\}$.  We find that $D_{wm}$ dominates $D_{rd}$ when $\min_{\gamma}P(\theta \in K)$ is substantially below $1 - \alpha$.  In other words, the main cause for $\min_{\gamma}P(\theta \in K)$ being substantially below $1 - \alpha$ is that the preliminary t test sometimes chooses the wrong model.

\bigskip

\noindent \textbf{5. Discussion}


\medskip

We have defined and shown how to assess $D_{wm}$,
the deficit in coverage probability due to the possible choice of the wrong model,
and $D_{rd}$, the deficit in coverage probability due to re-use of the data.
The numerical results presented in Figure 1 and the Supplementary Material show that the main cause for the coverage probability, minimized over $\gamma$, being substantially below nominal is that the preliminary t test sometimes chooses the wrong model.
This suggests that data-splitting will not help improve the coverage probability,  but that some form of
frequentist model averaging will help. The model averaged tail area confidence interval
of Turek and Fletcher (2012) is a promising method (Kabaila, Welsh and Abeysekera, 2016, and Kabaila, Welsh and Mainzer, 2016 and Kabaila, 2017).

\bigskip

\noindent \textbf{Acknowledgement}

This work was supported by an Australian Government Research Training Program Scholarship.

\bigskip

\noindent \textbf{Appendix: Proof of Theorem 1}

\medskip

Let 
\begin{equation*}
Z = \dfrac{G - \rho(H - \gamma)}{(1 - \rho^2)^{1/2}}.
\end{equation*}
It follows from 
\eqref{(G,H)dist}
that $Z$, $H$ and $W$ are independent random variables and that $Z \sim N(0, 1)$.

The event $\left\{ \theta \in J \right\} = \big\{ g_1(W, H/W) \leq Z \leq g_2(W, H/W)  \big\}$.
Thus
\begin{align*}
P(\theta \in J) &= P \Big(g_1(W, H/W) \leq Z \leq g_2( W, H/W) \Big) \\[1.2ex]
&= \int_0^{\infty} \int_{-\infty}^{\infty} P \Big( g_1(w, h/w) \leq Z \leq g_2(w, h/w) \Big) \, \phi(h - \gamma) \, dh \, f_W(w) \, dw \\[1.2ex]
&=\int_0^{\infty} \int_{-\infty}^{\infty} \Psi \big( w, h/w \big) \, \phi(h - \gamma) \, dh \, f_W(w) \, dw.
\end{align*}

The event $\big\{ |T| < t_{m, \, 1 - \widetilde{\alpha}/2} \big\} = \big\{ -t_{m, \, 1 - \widetilde{\alpha}/2} \, W \leq H \leq t_{m, \, 1 - \widetilde{\alpha}/2} \, W \big\}$.
Hence $P\big(\theta \in J, \, |T| < t_{m, \, 1 - \widetilde{\alpha}/2} \big)$ is equal to 
\begin{align*}
& P\Big(  g_1(W, H/W) \leq Z \leq g_2(W, H/W), \, -t_{m, \, 1 - \widetilde{\alpha}/2} \, W \leq H \leq t_{m, \, 1 - \widetilde{\alpha}/2} \, W  \Big) \\[1.2ex]
&= \int_0^{\infty} \int_{-t_{m, 1 - \widetilde{\alpha}/2} \, w}^{t_{m, 1 - \widetilde{\alpha}/2} \, w} P \Big( g_1(w, h/w) \leq Z \leq g_2(w, h/w) \Big) \, \phi \big( h - \gamma \big) \, dh \, f_W(w) \, dw.
\end{align*}
Changing the variable of integration of the inner integral to $y = h \big/ (t_{m, 1 - \widetilde{\alpha}/2} \, w)$, this expression becomes
\begin{equation*}
t_{m, \, 1 - \widetilde{\alpha}/2} \,
\int_0^{\infty} \int_{-1}^1 \Psi \big( w, t_{m, \, 1 - \widetilde{\alpha}/2} \, y \big)  \,  \phi \big(t_{m, \, 1- \widetilde{\alpha}/2} \, w \, y - \gamma \big) \, dy \,  w \, f_W(w)  \, dw.
\end{equation*}

The event $\{\theta \in I \} =
\big\{ -t_{m, \, 1 - \widetilde{\alpha}/2} \, W \leq G \leq t_{m, \, 1 - \widetilde{\alpha}/2} \, W \big\}$.
Therefore $P(\theta \in I, \, |T| < t_{m, \, 1 - \widetilde{\alpha}/2})$ is equal to 
\begin{equation*}
P \left( -t_{m, \, 1- \alpha/2} \leq \dfrac{G}{W} \leq t_{m, \, 1 - \alpha/2},  \, -t_{m, \, 1 - \widetilde{\alpha}/2} \leq \dfrac{H}{W} \leq t_{m, \, 1 - \widetilde{\alpha}/2} \right).
\end{equation*}
%


\baselineskip=18pt

\noindent \textbf{References}

\medskip

\rf Buckland, S.T., Burnham, K.P. and Augustin, N.H., 1997. Model selection: an integral
part of inference. \textsl{Biometrics}, 53, 603--618.

\medskip

\rf Chatfield, D., 1995.  Model uncertainty, data mining and statistical inference.  \textsl{Journal of the Royal Statistical Society, Series A}, 158, 419--466.

\medskip

\rf Efron, B., 2014. Estimation and accuracy after model selection. \textsl{Journal of the American Statistical Association}, 109, 991--1022.

\medskip

\rf Fletcher, D. and Turek, D., (2011). Model-averaged profile likelihood confidence intervals.
\textsl{Journal of Agricultural,  Biological and  Environmental Statistics}, 17, 38--51.

\medskip

\rf Hjort, N.L. and Claeskens, G., 2003. Frequentist model average estimators. 
\textsl{Journal of the American Statistical Association}, 98, 879--899.

\medskip

\rf Hurvich, C. M. and Tsai, C. L., 1990.  The impact of model selection on inference in linear regression.  \textsl{The American Statistician}, 44, 214--217.

\medskip 

\rf Kabaila, P., 2009.  The coverage properties of confidence regions after model selection.  \textsl{International Statistical Review}, 77, 405 -- 414.

\medskip 

\rf Kabaila, P., 2017.  On the minimum coverage probability of model averaged tail
area confidence intervals. Accepted for publication in \textsl{Canadian Journal of Statistics}.

\medskip

\rf Kabaila, P. and Giri, K., 2009.  Upper bounds on the minimum coverage probability of confidence intervals in regression after model selection.  \textsl{Australian \& New Zealand Journal of Statistics}, 51, 271--287.

\medskip

\rf Kabaila, P., Welsh, A.H. and Abeysekera, W. 2016. Model-averaged confidence intervals. 
\textsl{Scandinavian Journal of Statistics},
43, 35--48.

\medskip

\rf Kabaila, P., Welsh, A.H. and Mainzer, R., 2016. The performance of model averaged tail area confidence intervals. 
\textsl{Communications in Statistics - Theory and Methods}, 
DOI: 10.1080/03610926.2016.1242741

\medskip

\rf Leeb, H. and P\"otscher, B.,  2005. Model selection and inference: facts and fiction.  \textsl{Econometric Theory} 21, 21--59.

\medskip

\rf Miller, A., 2002. Subset Selection in Regression, 2nd ed.
Chapman \& Hall/CRC, Boca Raton FL.

\medskip

\rf Sheather, S. J., 2009.  A Modern Approach to Regression with R.  Springer, New York.

\medskip

\rf Turek, D. and Fletcher, D., 2012.  Model-averaged Wald intervals.  \textsl{Computational Statistics and Data Analysis}, 56, 2809--2815.

\end{document}